# ON APPROXIMATE PATTERN MATCHING FOR A CLASS OF GIBBS RANDOM FIELDS

By Jean-Rene Chazottes, Frank Redig and Evgeny Verbitskiy

*CNRS-Ecole Polytechnique, Leiden University and Philips Research*

We prove an exponential approximation for the law of approximate occurrence of typical patterns for a class of Gibssian sources on the lattice $\mathbb{Z}^d$, $d \geq 2$. From this result, we deduce a law of large numbers and a large deviation result for the waiting time of distorted patterns.

**1. Introduction.** In recent years there has been growing interest in a detailed probabilistic analysis of pattern matching and approximate pattern matching. For example, in information theory, motivation comes from studying performance of idealized Lempel–Ziv coding schemes. In mathematical biology one likes to have accurate estimates for the probability that two (e.g., DNA) sequences agree in a large interval with some error-percentage. There is also considerable interest in the analysis of occurrence of patterns in the multi-dimensional setting, for example, in the context of video-image compression [2], and more generally, lossy data compression [5, 6, 10].

In this paper we study the following problem. Fix a pattern $A_n$ in a cubic box of size $n$. Given a configuration $\sigma$ of a Gibbs random field, what is the size of the "observation window" in which we do not necessarily see exactly this pattern for the first time, but any pattern obtained by distortion of the fixed pattern $A_n$? By this, we mean a pattern which contains a fixed fraction $\varepsilon$ of spins different from those of $A_n$. We are interested in the behavior of the volume of this observation window, which we call "approximate hitting-time," when $n$ grows.

Our main result (Theorem 2.6) can be phrased as follows. The distribution of the approximate hitting-time, when properly normalized, gets closer and closer to an exponential law. The normalization is the product of a certain parameter $\Lambda_n$ and the probability of the set of distorted patterns $[A_n]^\varepsilon$.

In fact, we get a precise control of the error term which allows us to derive two corollaries for the "approximate waiting-time": given a configuration $\eta$ randomly chosen from an ergodic Gibbs random field, we increase the observation window in a random configuration $\sigma$ drawn from the given Gibbs random field until we see approximately the pattern $\eta_{C_n}$. The first corollary implies a law of large numbers allowing to get the rate-distortion function almost surely from this approximate waiting-time. The second corollary is related to large deviation bounds. While the law of large numbers for approximate waiting-times appears in [6] (under different conditions), the large deviation result is new. We emphasize that Theorem 2.6 is a new result in the context of approximate pattern-matching.

We briefly indicate the key ingredients needed to prove this exponential approximation. First, we assume that the Gibbs random field satisfies a certain strong mixing condition (nonuniformly $\varphi$-mixing condition). For instance, this property holds for all Markov random fields which satisfy the Dobrushin uniqueness condition. The second key ingredient is a result by Chi [4] allowing one to obtain the rate distortion function "à la Shannon–McMillan–Breiman." We take advantage of our previous work [1] in which we deal with "exact" hitting-times. The proof of the main result of the present work readily follows a large part of the proof in [1], but there is a crucial step which is different (second moment estimate). Moreover, one has to restrict to "good" patterns: if a pattern has "too much overlap" with its translates by vectors of size of order $n$, then one cannot hope to obtain an exponential distribution. These good patterns are shown to be typical in the sense that their measure approaches one exponentially fast as $n$ diverges (Proposition 2.7). When we have a random field distributed according to a Bernoulli measure, the goodness assumption on patterns can be removed. In this case, we prove (Theorem 2.8) that for any pattern Theorem 2.6 applies. Surprisingly, our proof involves the strong invariance principle for simple random walks. We have no idea how to provide a simpler proof.

*Outline of the paper.* In the next section we set notation and definitions and state our main theorems. In Section 3 we apply the exponential approximation of the previous section to approximate waiting times for which we obtain a.s. strong approximation and large deviations results. In Section 4 we state our proofs.

**2. Set-up and main results.** For the sake of simplicity, we consider a $\{0, 1\}$-valued random field on the lattice $\mathbb{Z}^d$, $d \geq 2$. Our results hold for any finite alphabet as well. Configurations are denoted $\eta, \sigma, \omega$ and collected in the set $\Omega = \{0, 1\}^{\mathbb{Z}^d}$. $\Omega$ is provided with the Borel $\sigma$-field, and for $V \subseteq \mathbb{Z}^d$, $\mathcal{F}_V$ denotes the $\sigma$-algebra generated by $\{\sigma_x : x \in V\}$.

For a finite subset $V \subseteq \mathbb{Z}^d$ and configurations $\sigma, \eta \in \Omega$, we denote by

$$\Delta(V, \eta, \sigma) = \sum_{x \in V} |\eta_x - \sigma_x| \tag{2.1}$$

the number of mismatches between $\sigma$ and $\eta$ in the volume $V$, that is, the Hamming distance between $\eta_V$ and $\sigma_V$.

We denote by $C_n$ the $n$-cube $[0,n]^d \cap \mathbb{Z}^d$. An $n$-pattern is a map $A_n : C_n \to \{0,1\}$. It is naturally associated to its cylinder $[A_n] = \{\sigma \in \Omega : \sigma_{C_n} = A_n\}$. For a pattern $A_n$ and $x \in \mathbb{Z}^d$, we denote by $\theta_{-x} A_n$ the pattern supported on $C_n + x$ defined by $A_n(y+x) = A_n(y)$ ($y \in C_n$).

We let $[A_n]^\varepsilon$ denote the set of configurations which $\varepsilon$-match with $A_n$:

$$[A_n]^\varepsilon = \{\omega \in \Omega : \Delta(C_n, \omega, A_n) \leq \varepsilon |C_n|\}. \tag{2.2}$$

The set of configurations $[A_n]^\varepsilon$ can also be viewed as a set of $n$-patterns. With a slight abuse of notation, we will use the same symbol for the set of configurations and the set of $n$-patterns which are restrictions of configurations in $[A_n]^\varepsilon$ to $C_n$.

DEFINITION 2.1. The *approximate hitting-time* of $[A_n]^\varepsilon$ in a configuration $\sigma$ is defined as

$$\mathbf{T}_{[A_n]^\varepsilon}(\sigma) = \min\{|C_k| : k > 0, \exists x \in \mathbb{Z}^d, C_n + x \subseteq C_k \text{ and } \theta_{-x}\sigma \in [A_n]^\varepsilon\}. \tag{2.3}$$

In words, given a configuration, the approximate hitting-time of the distorted pattern $[A_n]^\varepsilon$ is by definition the smallest volume of a $k$-cube ("observation window") such that there is some translate of an $n$-cube, contained in the observation window, which "hits" the $\varepsilon|C_n|$-ball (in the Hamming distance) around $A_n$.

For $\varepsilon = 0$ (exact matching time or occurrence time of a pattern), we obtained in [1] an exponential approximation for the law of $\mathbf{T}_{[A_n]^\varepsilon}$ under the hypotheses of nonuniform $\varphi$-mixing and Gibbsianness of the random field. We recall here this mixing assumption. For $m > 0$, define

$$\varphi(m) = \sup \frac{1}{|A_1|} |\Pr(E_{A_1} | E_{A_2}) - \Pr(E_{A_1})|, \tag{2.4}$$

where the supremum is taken over all finite subsets $A_1, A_2$ of $\mathbb{Z}^d$, with $d(A_1, A_2) \geq m$ [as usual, $d(A_1, A_2) := \inf\{d(x,y) : x \in A_1, y \in A_2\}$ and $d(x,y) := \|x-y\|_\infty = \max_{1 \leq i \leq d} |x_i - y_i|$] and $E_{A_i} \in \mathcal{F}_{A_i}$, with $\Pr(E_{A_2}) > 0$. Note that this $\varphi(m)$ differs from the usual $\varphi$-mixing function since we divide by the size of the dependence set of the event $E_{A_1}$.

DEFINITION 2.2. A random field is nonuniformly exponentially $\varphi$-mixing if there exist constants $C_1, C_2 > 0$ such that

$$\varphi(m) \leq C_1 e^{-C_2 m} \qquad \text{for all } m > 0. \tag{2.5}$$

A typical example of a Gibbs field satisfying this assumption is the 2d-Ising model above critical temperature. In general, it is satisfied in the so-called high-temperature regime of Dobrushin uniqueness. We refer the reader to [8, 9] for more details on this and on Gibbs measures in general.

An important property of Gibbs measures is the so-called "finite energy" property. This means that there is a continuous version of the conditional probability $\mathbb{P}(\sigma_0 = 0 | \sigma_{\mathbb{Z}^d \setminus \{0\}})$ such that

$$\delta < \mathbb{P}(\sigma_0 = 0 | \sigma_{\mathbb{Z}^d \setminus \{0\}}) < (1 - \delta), \tag{2.6}$$

where $\delta \in (0, \frac{1}{2})$ is independent of $\sigma$. This immediately implies the existence of $\kappa > 0$ such that, for all $V \subseteq \mathbb{Z}^d$, and all $\eta \in \Omega$,

$$\mathbb{P}(\{\sigma : \sigma_V = \eta_V\}) \leq e^{-\kappa |V|}. \tag{2.7}$$

We will use the following estimate:

LEMMA 2.3. *Under the assumption that $\mathbb{P}$ is a Gibbs measure, there exist $\varepsilon_c > 0$ and $K = K(\varepsilon_c) > 0$ such that, for any pattern $A_n$ and any $\varepsilon < \varepsilon_c$,*

$$\mathbb{P}([A_n]^\varepsilon) \leq e^{-Kn^d}.$$

PROOF. This is an immediate consequence of the estimate (2.7) and the estimate

$$|[A_n]^\varepsilon| \leq \sum_{k=0}^{\varepsilon n^d} \binom{n^d}{k} \leq e^{n^d I(\varepsilon)},$$

with $I(\varepsilon) \downarrow 0$ if $\varepsilon \downarrow 0$. $\square$

Contrary to the situation for exact matching, we will need an assumption on the patterns in order to obtain an exponential law. This can be compared with the condition of not being "badly self-repeating" needed to obtain the exponential law for return times in [1]. As we shall see, being a "good" pattern is a typical property.

DEFINITION 2.4. Given $0 < \alpha < 1, 0 \leq \varepsilon < 1$, we say that an *n-pattern* $A_n$ *is $(\varepsilon, \alpha)$-good* if the set $[A_n]^\varepsilon \cap \theta_x [A_n]^\varepsilon$ is empty for all $x \in \mathbb{Z}^d$ such that $|x| \leq \alpha n$. The set of all $(\varepsilon, \alpha)$-good patterns is denoted by $\mathcal{G}_n(\varepsilon, \alpha)$. By abuse of notation, we use the same symbol for the set of configurations $\omega$ such that $\omega_{C_n}$ is $(\varepsilon, \alpha)$-good.

For $\varepsilon = 0$ and $\alpha < 1/2$, $\mathcal{G}_n(\varepsilon, \alpha)$ coincides with the set of nonbadly self-repeating patterns in [1], Definition 5.1.

We shall need a result by Chi [4] on the rate distortion function. We recall briefly the definition of the rate distortion function and refer the reader to [3] for more information and background and to [6] for a discussion on lossy data compression. Given a stationary and ergodic measure $\mathbb{Q}$ and a stationary and ergodic Gibbs measure $\mathbb{P}$, the rate distortion function $\mathcal{R}(\mathbb{Q}, \mathbb{P}, \varepsilon)$ is defined as follows:

$$(2.8) \qquad \mathcal{R}(\mathbb{Q}, \mathbb{P}, \varepsilon) = \lim_{n \to \infty} \mathcal{R}_n(\mathbb{Q}, \mathbb{P}, \varepsilon),$$

$$(2.9) \qquad \mathcal{R}_n(\mathbb{Q}, \mathbb{P}, \varepsilon) = \inf_{\mathbb{J}_n} \frac{1}{|C_n|} H(\mathbb{J}_n \parallel \mathbb{Q}_n \times \mathbb{P}_n),$$

where the infimum taken over all joint distributions $\mathbb{J}_n$ on $\{0,1\}^{n^d} \times \{0,1\}^{n^d}$ such that the $\{0,1\}^{n^d}$-marginal of $\mathbb{J}_n$ is $\mathbb{Q}_n$ and

$$\int \frac{\Delta(C_n, \omega, \sigma)}{|C_n|} \, d\mathbb{J}_n(\omega, \sigma) \leq \varepsilon.$$

$H(\mathbb{J}_n \parallel \mathbb{Q}_n \times \mathbb{P}_n)$ is the relative entropy between $\mathbb{J}_n$ and $\mathbb{Q}_n \times \mathbb{P}_n$.

We have the following key result which follows from [4] and [6], Theorem 25.

PROPOSITION 2.5. *Let $\mathbb{Q}$ be a stationary and ergodic measure and $\mathbb{P}$ be a stationary and ergodic Gibbs measure. Then*

$$(2.10) \quad \mathcal{R}(\mathbb{Q}, \mathbb{P}, \varepsilon) = -\lim_{n \to \infty} \frac{1}{|C_n|} \log \mathbb{P}([\omega_{C_n}]^\varepsilon) \qquad \mathbb{Q}\text{-almost-surely}.$$

*Moreover, $\mathcal{R}$ is a convex (and, hence, continuous) function of $\varepsilon$ and is nonzero in some interval $[0, \varepsilon_0)$.*

The property (2.10) is called the generalized asymptotic equipartition property in [6]. Throughout we will simply write $\mathcal{R}(\varepsilon)$ instead of $\mathcal{R}(\mathbb{Q}, \mathbb{P}, \varepsilon)$.

We can now state our main result.

THEOREM 2.6. *Suppose that $\mathbb{P}$ is a nonuniformly exponentially $\varphi$-mixing Gibbs measure and $\mathbb{Q}$ is a stationary and ergodic Gibbs measure. Assume that the rate distortion function (2.8) is strictly positive in $[0, \varepsilon_0)$. Then for all $\alpha \in (0, 1)$ and $\varepsilon > 0$ small enough, namely,*

$$\frac{\varepsilon}{\alpha} < \varepsilon_0,$$

*there exist $\Lambda_1, \Lambda_2, C, c \in (0, \infty)$, such that and for every $t > 0$, $n \geq 1$, and $\mathbb{Q}$-almost all $\omega$ with $\omega_{C_n} \in \mathcal{G}_n(\varepsilon, \alpha)$, the following estimate holds:*

$$(2.11) \qquad \left| \mathbb{P}\left( \mathbf{T}_{[\omega_{C_n}]^\varepsilon} > \frac{t}{\Lambda_n \mathbb{P}([\omega_{C_n}]^\varepsilon)} \right) - e^{-t} \right| \leq C e^{-ct} \, e^{-Kn^d},$$

*where* $\Lambda_n = \Lambda(\omega_{C_n})$ *is such that*

(2.12) $$\Lambda_1 \leq \Lambda_n \leq \Lambda_2.$$

Dependence of the parameters in Theorem 2.6 on $\varepsilon$ and $\alpha$ will be discussed after the proof; see Remark 4.1.

Let us briefly comment on the difference between Theorem 2.6 and the one obtained in [1] for exact matching, that is, the case corresponding to $\varepsilon = 0$. First of all, we need to restrict ourselves to special patterns, that is, $(\varepsilon, \alpha)$-good patterns, whereas in [1] result applies to all patterns. Second, the error term that we obtain in [1] is of the form $Ce^{-ct}\mathbb{P}([\omega_{C_n}])^\rho$, where $\rho > 0$. Of course, the factor $\mathbb{P}([\omega_{C_n}])^\rho$ is uniformly exponentially small for Gibbs measures. This is no longer true for $\mathbb{P}([\omega_{C_n}]^\varepsilon)$ if $\varepsilon$ is too large. This is precisely why we need Lemma 2.3. Third, a crucial step in the proof of Theorem 2.6, which differs slightly from that in [1] for the case $\varepsilon = 0$, involves Proposition 2.5. This explains why we need to restrict to typical configurations in the sense of this result.

Let us close this set of remarks by noticing that $\mathbb{Q}$ has to be a stationary and ergodic measure, but not necessarily Gibbsian. But for later use of Theorem 2.6, we shall also need the latter assumption, so we already impose it to state the theorem.

The following proposition shows that "$\omega_{C_n} \in \mathcal{G}(\varepsilon, \alpha)$," that is, that a pattern being $(\varepsilon, \alpha)$-good, is a typical property.

PROPOSITION 2.7. *Let $\mathbb{Q}$ be a stationary Gibbs measure. Then, if $\alpha < 1/2$ and $\varepsilon > 0$ is small enough, there exists $\nu > 0$ such that, for all $n \geq 1$,*

(2.13) $$\mathbb{Q}(\mathcal{G}_n(\varepsilon, \alpha)) > 1 - e^{-\nu n^d}.$$

It turns out that if the random field has a nontrivial dependence structure, then the restriction to $(\varepsilon, \alpha)$-good patterns is unavoidable. However, in the case of a random field distributed according to a Bernoulli measure, the exponential law (2.11) holds for *all* approximate patterns. This is expressed by the following theorem.

THEOREM 2.8. *If $\mathbb{P}$ is the Bernoulli measure with $\mathbb{P}(\sigma_0 = 1) = 1/2$, then (2.11) holds without the restriction that $\omega_{C_n}$ is $(\varepsilon, \alpha)$-good.*

**3. Approximate waiting-time fluctuations.** The purpose of this section is to derive two consequences of Theorem 2.6 and Proposition 2.7. The first one implies a strong law of large numbers for the approximate waiting-time. It was previously derived in [6] directly using the mixing property (2.5). The second one concerns large deviations of the approximate waiting-time and it is a new result. Given two configurations $\omega, \sigma$, the *approximate waiting-time* is $\mathbf{W}_n^\varepsilon(\omega, \sigma) := \mathbf{T}_{[\omega_{C_n}]^\varepsilon}(\sigma)$.

PROPOSITION 3.1. *Under the assumptions of Theorem* 2.6 *and Proposition* 2.7, *there exists* $\gamma_0 > 0$ *such that, for all* $\gamma > \gamma_0$,

$$(3.1) \qquad -\gamma \log n \leq \log(\mathbf{W}_n^\varepsilon(\omega,\sigma)\mathbb{P}([\omega_{C_n}]^\varepsilon)) \leq \log(\log n^\gamma)$$

$\mathbb{Q} \times \mathbb{P}$-*eventually almost surely. In particular,*

$$\lim_{n\to\infty} \frac{1}{|C_n|} \log \mathbf{W}_n^\varepsilon(\omega,\sigma) = \mathcal{R}(\mathbb{Q},\mathbb{P},\varepsilon), \qquad \mathbb{Q} \times \mathbb{P}\text{-almost surely.}$$

With Proposition 3.1, we recover the results of Theorems 26 and 27 in [6]. However, there is a substantial difference in conditions on random fields. We have to restrict ourselves to measures $\mathbb{Q}$ which are stationary and ergodic Gibbs measures, while in [6] $\mathbb{Q}$ is only assumed to be stationary and ergodic. On the other hand, we permit $\mathbb{P}$ to be Gibbsian, while in [6] $\mathbb{P}$ must be a Bernoulli measure. The reason for our assumptions on $\mathbb{Q}$ is that Proposition 2.7 is valid for Gibbs measures. We do not know if it can be extended to more general situations.

Let us also remark that, by a basic result in Probability Theory, this strong approximation implies that if a central limit theorem holds for $-1/|C_n| \times \log \mathbb{P}([\omega_{C_n}]^\varepsilon)$, then it holds also for $(1/|C_n|) \log \mathbf{W}_n^\varepsilon(\omega,\sigma)$. Unfortunately, the former seems to be a difficult issue, except in the i.i.d. case. We refer the reader to [6] for some results in that direction.

We have the following (partial) large deviation results. We first need the following lemma showing that we can define the generalized conditional $q$-order Rényi entropy for Gibbs random fields. This was first done in [11] for ($\alpha$-mixing) stochastic processes ($d=1$) with the difference that here we need to condition on ($\varepsilon,\alpha$)-good patterns and use the Gibbs property instead of mixing.

LEMMA 3.2. *Let* $\mathbb{Q}, \mathbb{P}$ *be stationary Gibbs measures and assume that* $\alpha < 1/2$ *and* $0 \leq \varepsilon < 1$. *Then, for all* $q \in \mathbb{R}$, *the following function is well defined:*

$$(3.2) \quad \mathcal{E}_\varepsilon(q) := \mathcal{E}_\varepsilon(q;\mathbb{Q},\mathbb{P}) = \lim_{n\to\infty} \frac{1}{|C_n|} \log \int \mathbb{P}([\omega_{C_n}]^\varepsilon)^q \, d\mathbb{Q}_{\mathcal{G}_n(\varepsilon,\alpha)}(\omega).$$

($\mathbb{Q}_{\mathcal{G}_n(\varepsilon,\alpha)}$ *denotes the measure* $\mathbb{Q}$ *conditioned on the set of good patterns.*)

The generalized $q$-order Rényi entropy should be defined as $-\mathcal{E}_\varepsilon(-q)/q$. We now have the following theorem. By $a_n \approx b_n$, we mean that

$$\max\{a_n/b_n, \ b_n/a_n\}$$

is bounded from above.

THEOREM 3.3. *Let $\mathbb{P}$ be a nonuniformly exponentially $\varphi$-mixing Gibbs measure and $\mathbb{Q}$ a stationary and ergodic Gibbs measure. If $\varepsilon > 0$ is small enough, then for any $\alpha_0 \leq \alpha < 1/2$, we have*

$$\iint (\mathbf{W}_n^\varepsilon(\omega, \sigma))^q \, d\mathbb{Q}_{\mathcal{G}_n(\varepsilon,\alpha)}(\omega) \, d\mathbb{P}(\sigma)$$
(3.3)
$$\approx \int \mathbb{P}([\omega_{C_n}]^\varepsilon)^{-q} \, d\mathbb{Q}_{\mathcal{G}_n(\varepsilon,\alpha)}(\omega) \quad \text{if } q \geq -1$$

*and*

$$\iint (\mathbf{W}_n^\varepsilon(\omega, \sigma))^q \, d\mathbb{Q}_{\mathcal{G}_n(\varepsilon,\alpha)}(\omega) \, d\mathbb{P}(\sigma)$$
(3.4)
$$\approx \int \mathbb{P}([\omega_{C_n}]^\varepsilon) \, d\mathbb{Q}_{\mathcal{G}_n(\varepsilon,\alpha)}(\omega) \quad \text{if } q < -1.$$

*In particular,*

(3.5)
$$\lim_{n\to\infty} \frac{1}{|C_n|} \log \iint (\mathbf{W}_n^\varepsilon(\omega, \sigma))^q \, d\mathbb{Q}_{\mathcal{G}_n(\varepsilon,\alpha)}(\omega) \, d\mathbb{P}(\sigma)$$
$$= \begin{cases} \mathcal{E}_\varepsilon(-q), & \text{if } q \geq -1, \\ \mathcal{E}_\varepsilon(1), & \text{if } q < -1. \end{cases}$$

It follows from this theorem that Theorem 4.5.20 in [7] applies to $\{1/|C_n| \times \log \mathbf{W}_n^\varepsilon(\eta, \sigma)\}_n$. However, to obtain a full large deviation principle, we would need to know under which conditions the function $q \mapsto \mathcal{E}_\varepsilon(q)$ is, for instance, differentiable for $q > -1$ (and for $\varepsilon$ small enough). If that were the case, we would have a large deviation principle with a rate function given by the Legendre transform of $\mathcal{E}_\varepsilon(-q)$.

## 4. Proofs.

4.1. *Proof of Theorem* 2.6. The proof of Theorem 2.6 is quite similar to the proof of exponential law in [1]. We describe briefly the common approach and indicate the differences. We also provide the necessary modifications of the proof.

It is well known that a random variable $Z$ has an exponential distribution if and only if

$$\mathbb{P}(Z > s + t | Z > t) = \mathbb{P}(Z > s)$$

or, equivalently,

$$\mathbb{P}(Z > s + t) = \mathbb{P}(Z > s)\mathbb{P}(Z > t).$$

The basic ingredient of the proof in [1] was Lemma 4.4 ("Iteration Lemma"). This result establishes that, for a pattern $A_n$ and any finite number of cubes $C_i \subseteq \mathbb{Z}^d$, $i = 1, \ldots, k$, with equal volumes

$$|C_i| = \left(\frac{1}{\mathbb{P}(A_n)}\right)^\gamma,$$

we have

(4.1) $\mathbb{P}\left(A_n \text{ does not occur in } \bigcup_{i=1}^k C_i\right) \approx \mathbb{P}(A_n \text{ does not occur in } C_1)^k.$

In [1] we also observed that the Iteration Lemma remains valid if a pattern $A_n$ is replaced by the event $[A_n]^\varepsilon$, with $[A_n]^\varepsilon$ not occuring in volume $V$ if any pattern $B_n \in [A_n]^\varepsilon$ does not occur in volume $V$.

Another important ingredient of the proof is the control of the parameter of the exponential distribution. Lemma 4.3 ("The parameter") in [1] concerns nontriviality of the parameter $\Lambda_n$, that is, the fact that it is neither null nor infinite. To prove Lemma 4.3, we established a *uniform* second moment estimate for the number of occurrences of a pattern $A_n$ in a configuration $\sigma$ restricted to a box that has later to be taken of size $1/\mathbb{P}([A_n])$. It is the proof of this second moment estimate that we have to modify completely. In Remark 4.1 in [1], we noticed that if $E_n \in \mathcal{F}_{C_n}$ are events such that $\mathbb{P}(E_n) < e^{-cn^d}$ for some $c > 0$, and such that

(4.2) $$\limsup_{n \to \infty} \sum_{0 < |x| < n} \frac{\mathbb{P}(E_n \cap \theta_x E_n)}{\mathbb{P}(E_n)} < \infty,$$

then this implies, together with the mixing property (2.5) and the Gibbs property (2.7), that the desired uniform second moment estimate holds. In turn, this implies the nontriviality of the parameter (2.12) (Lemma 4.3 in [1]).

Thus, we turn to prove (4.2) when the event $E_n$ is $[A_n]^\varepsilon$, where $A_n$ is a good and typical pattern. We assume that patterns $A_n$ are such that $\bigcap_n A_n = \{\sigma\}$, with $\sigma$ chosen in the set with $\mathbb{Q}$-measure one from Proposition 2.5, and such that $A_n$ is good in the sense of Definition 2.4.

We have to show for patterns $A_n \in \mathcal{G}_n(\varepsilon, \alpha)$ with $\varepsilon/\alpha < \varepsilon_0$ that there exists a finite number $C(\varepsilon, \alpha)$ such that, for all $n$,

(4.3) $$\sum_{0 < |x| \leq n} \frac{\mathbb{P}([A_n]^\varepsilon \cap \theta_x[A_n]^\varepsilon)}{\mathbb{P}([A_n]^\varepsilon)} \leq C(\varepsilon, \alpha).$$

First of all, since $A_n \in \mathcal{G}_n(\varepsilon, \alpha)$ (see Definition 2.4), the terms corresponding to $x$ with $|x| < \alpha n$ are equal to 0. Therefore, we have to estimate the

sum

(4.4) $$\sum_{\alpha n \leq |x| \leq n} \frac{\mathbb{P}([A_n]^\varepsilon \cap \theta_x[A_n]^\varepsilon)}{\mathbb{P}([A_n]^\varepsilon)}.$$

Note that, for $x$ with $|x| \geq \alpha n$, the intersection $(C_n + x) \cap C_n$ is not very large:

$$|(C_n + x) \cap C_n| \leq (1-\alpha)n^d.$$

Note also that $\Delta(V, \omega, A_n)$ denotes the number of differences between $\omega$ and $A_n$ in the volume $V$, see (2.1). Then we can write

(4.5) $$\mathbb{P}([A_n]^\varepsilon \cap \theta_x[A_n]^\varepsilon)$$
$$= \mathbb{P}(\omega : \Delta(C_n, \omega, A_n) \leq \varepsilon n^d \cap \Delta(C_n + x, \omega, \theta_{-x}A_n) \leq \varepsilon n^d)$$

where, by $\theta_{-x}A_n$, we mean $\theta_{-x}A_n(y + x) = A_n(y)$, $y \in C_n$. For the sake of convenience, we simply write $C$ for $C_n$ and $C_x$ for $C_n + x$ in the course of this proof. We also introduce the short-hand notation

(4.6) $$\begin{aligned} S_1 &= \Delta(C \setminus C_x, \omega, A_n), \\ S_2 &= \Delta(C \cap C_x, \omega, A_n), \\ S_3 &= \Delta(C \cap C_x, \omega, \theta_{-x}A_n), \\ S_4 &= \Delta(C \setminus C_x, \omega, \theta_{-x}A_n). \end{aligned}$$

With this notation what we have to estimate is

(4.7) $$\sum_{\alpha n \leq |x| \leq n} \frac{\mathbb{P}([A_n]^\varepsilon \cap \theta_x[A_n]^\varepsilon)}{\mathbb{P}([A_n]^\varepsilon)}$$
$$= \sum_{\alpha n \leq |x| \leq n} \mathbb{P}(S_3 + S_4 \leq \varepsilon n^d | S_1 + S_2 \leq \varepsilon n^d).$$

The following estimate is a corollary of [4] and a basic property of a Gibbs measures: for any configuration $\xi$,

(4.8) $$\mathbb{P}(\{\omega : \Delta(V_n, \omega, \sigma) \leq \varepsilon |V_n|\} | \xi_{V_n^c}) \leq \exp(-|V_n| \, \mathcal{R}(\varepsilon) + c|\partial V_n|).$$

Indeed, the unconditioned statement is proved in [4], and conditioning can at most introduce a term of order $\exp(c|\partial V_n|)$.

We proceed as follows:

$$\mathbb{P}((S_1 + S_2) \leq \varepsilon n^d \cap (S_3 + S_4) \leq \varepsilon n^d)$$
$$\leq \mathbb{P}((S_1 + S_2) \leq \varepsilon n^d \cap S_4 \leq \varepsilon n^d)$$
$$\leq \sup_\xi \mathbb{P}(S_4 \leq \varepsilon n^d | \xi_{\mathbb{Z}^d \setminus (C_x \setminus C)}) \mathbb{P}([A_n]^\varepsilon)$$
$$\leq \exp\left(-\alpha n^d \mathcal{R}\left(\frac{\varepsilon}{\alpha}\right) + cn^{d-1}\right) \mathbb{P}([A_n]^\varepsilon).$$

Therefore,

$$\sum_{\alpha n \leq |x| \leq n} \frac{\mathbb{P}([A_n]^\varepsilon \cap \theta_x [A_n]^\varepsilon)}{\mathbb{P}([A_n]^\varepsilon)} \leq n^d \exp\left(-\alpha n^d \mathcal{R}\left(\frac{\varepsilon}{\alpha}\right) + cn^{d-1}\right) =: C_n(\varepsilon, \alpha).$$

Taking into account that $\varepsilon/\alpha < \varepsilon_0$, and, hence, $\mathcal{R}(\varepsilon/\alpha) > 0$, we conclude that $C_n(\varepsilon, \alpha) \to 0$ as $n \to \infty$, and, hence,

$$C(\varepsilon, \alpha) = \sup_n C_n(\varepsilon, \alpha)$$

is finite. This completes the proof.

REMARK 4.1. The parameters of Theorem 2.6 depend on the choice of $\varepsilon$ and $\alpha$. The most interesting is the dependence of $\Lambda_1$ and $\Lambda_2$. Lemma 4.3 in [1] in fact shows that a uniform choice $\Lambda_2 = 2$ suffices. A more interesting question is whether we can give a uniform bound on $\Lambda_1$ for a large set of $\varepsilon$ and $\alpha$. The present modification of the second moment estimate, together with the rest of Lemma 4.3 in [1], which remains unchanged, gives that, for some $c$, dependent on $\varepsilon$ alone, the following choice of $\Lambda_1 = \Lambda_1(\varepsilon, \alpha)$ will suffice:

$$\Lambda_1 = \frac{1}{c + C(\varepsilon, \alpha)}.$$

The rate distortion function $\mathcal{R}$ is a monotonically decreasing function. Hence, for a fixed $\varepsilon > 0$, $\alpha \mathcal{R}(\frac{\varepsilon}{\alpha})$ is a monotonically increasing function of $\alpha$, and finally, $C(\varepsilon, \alpha)$ is monotonically decreasing in $\alpha$. Therefore, if $\varepsilon < \varepsilon_0$, then for all $\alpha > \alpha_0 := 0.99 \frac{\varepsilon}{\varepsilon_0}$,

$$\Lambda_1(\varepsilon, \alpha) \geq \Lambda_1(\varepsilon, \alpha_0) > 0.$$

Therefore, for a fixed $\varepsilon > 0$, we obtain a uniform (in $\alpha$) bound on the parameter $\Lambda_1$.

4.2. *Proof of Proposition* 2.7. For $\varepsilon = 0$, we know that most patterns are $(0, \alpha)$-good for any $\alpha < 1$. Indeed, it is proved in [1] (Lemma 5.3) that $\mathbb{Q}(\mathcal{G}_n(\varepsilon, \alpha)) \geq 1 - e^{-\kappa' n^d}$, for some $\kappa' > 0$.

Let us now argue that for small $\varepsilon$ this is still the case. Suppose $\alpha < 1/2$, that is, we are going to consider vectors $x \in \mathbb{Z}^d$ such that $|x| \leq \frac{n}{2}$. An element $A$ of $[A_n]^\varepsilon \cap \theta_x [A_n]^\varepsilon$ satisfies

(4.9) $$\sum_{y \in C_x \cap C} |A(y) - A(y-x)| \leq 2\varepsilon n^d.$$

(Recall that $C = C_n$ and $C_x = C_n + x$.) This implies that there exists a set $V_n \subseteq C$ and a disjoint translate $V_n + z \subseteq C$ such that $|V_n| > (1/2)^d n^d$ such that $\theta_{-z} A_{V_n+z}$ matches with error fraction $2^{d+1}\varepsilon$ with $A_{V_n}$; this can be made as small as $e^{-\nu n^d}$, for some $\nu > 0$, for $\varepsilon$ sufficiently small uniformly in $A_{V_n}$ by Lemma 2.3. Therefore, we obtain that

$$\mathbb{Q}(\mathcal{G}(\varepsilon, \alpha_0)) > 1 - e^{-\nu n^d} \tag{4.10}$$

for all $\alpha < 1/2$ and $\varepsilon$ small enough.

4.3. *Proof of Theorem* 2.8. We consider the case $d = 1$ only, because the case $d \geq 2$ is completely analogous. Start with the particular pattern $A_n = 0 \cdots 0$ that we simply denote by $0_n$. The difficulty with this "bad pattern" comes from the fact that the second moment estimate does not apply, because (4.2) fails. Therefore, we have to prove by other means that there exists $\delta > 0$ such that, for all $n \in \mathbb{N}$,

$$\delta < \mathbb{P}\bigg(\mathbf{T}_{[0_n]^{\varepsilon}} > \frac{1}{\mathbb{P}([0_n]^{\varepsilon})}\bigg) < 1 - \delta, \tag{4.11}$$

which would imply the nontriviality of the parameter $\Lambda_n$. We will first show that there exists a sequence $k_n \uparrow \infty$ such that

$$\delta < \mathbb{P}(\mathbf{T}_{[0_n]^{\varepsilon}} > k_n) < 1 - \delta. \tag{4.12}$$

It will then follow easily from the Bernoulli character of $\mathbb{P}$ that $k_n$ does not depend on the choice of the pattern, that is, (4.12) holds with the same $k_n$ for any pattern $A_n$. Then we can apply Theorem 2.6 for good patterns, and obtain $k_n = 1/\mathbb{P}([A_n]^{\varepsilon}) = 1/\mathbb{P}([0_n]^{\varepsilon})$. We have the following identities:

$$\begin{aligned}
\mathbb{P}(\mathbf{T}_{[0_n]^{\varepsilon}} \leq k_n) &= \mathbb{P}\bigg(\min_{k=0}^{k_n} \sum_{i=k}^{k+n} \omega_i \leq n\varepsilon\bigg) \\
&= \mathbb{P}\bigg(\max_{k=0}^{k_n} \sum_{i=k}^{k+n} (1 - 2\omega_i) \geq (1 - 2\varepsilon)n\bigg) \\
&= \mathbb{P}\bigg(\max_{k=0}^{k_n} (S_{k+n} - S_k) \geq (1 - 2\varepsilon)n\bigg),
\end{aligned} \tag{4.13}$$

where $S_n$ is the position of a simple random walk on $\mathbb{Z}$ (with $S_0 = 0$) after $n$ steps. By Theorem 7.23 in [12], together with the strong invariance principle ([12], page 53), we have

$$\max_{k=0}^{k_n} (S_{k+n} - S_k) = a \log k_n + b \log \log k_n + c + o(1) + X, \tag{4.14}$$

where $X$ is a random variable with a Gumbel distribution. Therefore, if we choose $k_n$ such that

$$(1 - 2\varepsilon)n = a \log k_n + b \log \log k_n + c + o(1), \tag{4.15}$$

then (4.12) holds.

If we now choose any other pattern $A_n$, then, under $\mathbb{P}$, $S_n = 2 \sum_{i=0}^{n} (1/2 - \sigma_i - A_n(i))$ is again distributed as a simple random walk, so we find the same $k_n$, which completes the proof of the theorem.

4.4. *Proof of Proposition* 3.1. By using Theorem 2.6, we immediately get

$$\mathbb{Q} \times \mathbb{P}\{(\omega,\sigma) : \log(\mathbf{W}_n^\varepsilon(\omega,\sigma)\mathbb{P}([\omega_{C_n}]^\varepsilon)) > \log t\}$$
$$= \int d\mathbb{Q}(\omega)\ \mathbb{P}\{\sigma : \log(\mathbf{T}_{[\omega_{C_n}]^\varepsilon}(\sigma)\mathbb{P}([\omega_{C_n}]^\varepsilon)) > \log t\}$$
$$\leq e^{-\Lambda_1 t} + Ce^{-Kn^d} + \sum_{A_n \in \mathcal{G}_n^c(\varepsilon,\alpha)} \mathbb{Q}([A_n]).$$

Now we choose $t = t_n = \log(n^\gamma)$ with $\gamma > 0$ such that $\Lambda_1 \gamma > 1$. This makes the first term in the right–hand side summable in $n$. The last one equals the $\mathbb{Q}$-measure of the complement of $\mathcal{G}_n(\varepsilon, \alpha)$, which is less than $e^{-\nu n^d}$ by Proposition 2.7. We thus get the upper bound in (3.1) by an application of the Borel–Cantelli lemma.

Now we turn to prove the lower bound in (3.1). Proceeding as before, we get

$$\mathbb{Q} \times \mathbb{P}\{(\omega,\sigma) : \log(\mathbf{W}_n^\varepsilon(\omega,\sigma)\mathbb{P}([\omega_{C_n}]^\varepsilon)) \leq \log t\}$$
$$\leq 1 - e^{-\Lambda_2 t} + Ce^{-Kn^d} + \sum_{A_n \in \mathcal{G}_n^c(\varepsilon,\alpha)} \mathbb{Q}([A_n])$$
$$\leq \Lambda_2 t + Ce^{-Kn^d} + e^{-\nu n^d}.$$

We have used Theorem 2.6 and Proposition 2.7. We now choose $t = t_n = n^{-\gamma}$, with $\gamma > 1$, to get a summable upper bound in $n$ for the above probability. An application of Borel–Cantelli lemma gives the desired result and the proof of the proposition is complete.

4.5. *Proof of Lemma* 3.2. We only consider the case $q > 0$ leaving the (very similar) proof for the case $q < 0$ to the reader. Let $\mathcal{S}_\square$ be the system of all rectangular boxes of the form

$$V = \mathbb{Z}^d \cap \prod_{k=1}^d [m_k, n_k] \qquad \text{with } m_k, n_k \in \mathbb{Z},\ m_k \leq n_k.$$

Before proceeding, we have to extend Definition 2.4 somewhat. We will denote by $\mathcal{G}_V(\varepsilon, \alpha)$ the set of good patterns supported on $V \in \mathcal{S}_\square$. We shall need Proposition 2.7, which remains valid if one replaces $\mathcal{G}_n(\varepsilon, \alpha)$ with $\mathcal{G}_V(\varepsilon, \alpha)$ and $n$ by $|V|$ in (2.13).

We are going to prove that the function $a : \mathcal{S}_\square \to (-\infty, +\infty)$ defined as

$$a(V) := -\log \int \mathbb{P}^q([\sigma_V]^\varepsilon)\, d\mathbb{Q}_{\mathcal{G}_{V \cup V'}(\varepsilon,\alpha)}(\sigma)$$

satisfies the following approximate sub-additive property:

$$a(V \cup V') \leq a(V) + a(V') + C\,|\partial(V \cup V')|$$

for all $V, V' \in \mathcal{S}_\square$ such that $V \cup V' \in \mathcal{S}_\square$ and $V \cap V' = \varnothing$, where $C$ is a constant (depending on $q$), and where $\partial V$ denotes the boundary of $V$. Of course, $|\partial(V \cup V')|$ is a surface order correction. If such a property holds [together with $a(V + x) = a(V)$, for all $x \in \mathbb{Z}^d$, $V \in \mathcal{S}_\square$ which is obvious by stationarity of the measure], then a generalized sub-additive lemma, obtained as a combination of a lemma found in [8] and another one given in [7], will guarantee that

$$\lim_{n \to \infty} \frac{a(C_n)}{|C_n|}$$

exists, as we wish. For all $q \in \mathbb{R}$, $V, V' \in \mathcal{S}_\square$ such that $V \cup V' \in \mathcal{S}_\square$ and $V \cap V' = \varnothing$, we have the following:

$$\mathbb{P}^q([\sigma_{V \cup V'}]^\varepsilon) = \left( \sum_{\omega_{V \cup V'} \in [\sigma_{V \cup V'}]^\varepsilon} \mathbb{P}([\omega_{V \cup V'}]) \right)^q$$

$$\geq e^{K_1 |\partial(V \cup V')|} \left( \sum_{\omega_{V \cup V'} \in [\sigma_{V \cup V'}]^\varepsilon} \mathbb{P}([\omega_V]) \mathbb{P}([\omega_{V'}]) \right)^q$$

$$\geq e^{K_2 |\partial(V \cup V')|} \left( \sum_{\omega_V \in [\sigma_V]^\varepsilon} \mathbb{P}([\omega_V]) \right)^q \left( \sum_{\omega_{V'} \in [\omega_{V'}]^\varepsilon} \mathbb{P}([\omega_{V'}]) \right)^q$$

$$= e^{K_2 |\partial(V \cup V')|} \mathbb{P}^q([\sigma_{V \cup V'}]^\varepsilon) \mathbb{P}^q([\sigma_{V \cup V'}]^\varepsilon),$$

where $K_1, K_2$ are constants. The first inequality follows from the Gibbs property and the second one is a simple consequence of the Hamming distance property. To complete the proof, we again use the Gibbs property to get

$$\int \mathbb{P}^q([\sigma_{V \cup V'}]^\varepsilon) \, d\mathbb{Q}_{\mathcal{G}_{V \cup V'}(\varepsilon, \alpha)}(\sigma)$$

$$= \sum_{\omega_{V \cup V'} \in \{0,1\}^{V \cup V'}} \mathbb{P}^q([\omega_{V \cup V'}]^\varepsilon) \mathbb{Q}_{\mathcal{G}_{V \cup V'}(\varepsilon, \alpha)}([\omega_{V \cup V'}])$$

$$\geq \mathbb{Q}(\mathcal{G}_{V \cup V'}(\varepsilon, \alpha)) e^{K_3 |\partial(V \cup V')|}$$

$$\times \sum_{\omega_V \in \{0,1\}^V} \mathbb{P}^q([\omega_V]^\varepsilon) \mathbb{Q}_{\mathcal{G}_{V \cup V'}(\varepsilon, \alpha)}([\omega_V])$$

$$\times \sum_{\omega_{V'} \in \{0,1\}^{V'}} \mathbb{P}^q([\omega_{V'}]^\varepsilon) \mathbb{Q}_{\mathcal{G}_{V \cup V'}(\varepsilon, \alpha)}([\omega_{V'}])$$

$$\geq \tfrac{1}{2} e^{K_3 |\partial(V \cup V')|} \int \mathbb{P}^q([\sigma_V]^\varepsilon) \, d\mathbb{Q}_{\mathcal{G}_V(\varepsilon, \alpha)}(\sigma)$$

$$\times \int \mathbb{P}^q([\sigma_{V'}]^\varepsilon) \, d\mathbb{Q}_{\mathcal{G}_{V'}(\varepsilon, \alpha)}(\sigma),$$

where $K_3$ is a constant. The second inequality is the consequence of Proposition 2.7 if $|V \cup V'|$ is large enough. The lemma is proved.

4.6. *Proof of Theorem* 3.3. Since the proof of this theorem is very similar to that of Theorem 2.7 in [1], we only sketch it to indicate the little differences between them.

The starting point is of course to write

$$\iint (\mathbf{W}_n^\varepsilon(\omega,\sigma))^q \, d\mathbb{Q}_{\mathcal{G}_n(\varepsilon,\alpha)}(\omega) \, d\mathbb{P}(\sigma) = \int d\mathbb{Q}_{\mathcal{G}_n(\varepsilon,\alpha)}(\omega) \, \int \mathbf{T}_{[\omega_{C_n}]^\varepsilon}^q(\sigma) \, d\mathbb{P}(\sigma).$$

Then we can mimic the proof of Theorem 2.7 in [1] by using Theorem 2.6 and the analog of Lemma 4.3 in [1], which holds true when $\mathbf{T}_{[\omega_{C_n}]}$ is replaced by $\mathbf{T}_{[\omega_{C_n}]^\varepsilon}$, provided that $\omega_{C_n}$ be an $(\varepsilon,\alpha)$-good pattern (see the beginning of the proof of Theorem 2.6), and $\omega$ be $\mathbb{Q}$-typical in the sense of Proposition 2.5. Notice that we integrate with respect to the conditional measure $\mathbb{Q}_{\mathcal{G}_n(\varepsilon,\alpha)}$ which takes care of these two properties.

**Acknowledgment.** We thank Z. Chi for providing us with his preprint [4].

J.-R. Chazottes  
CPhT, CNRS-Ecole Polytechnique  
91128 Palaiseau Cedex  
France  
E-mail: jeanrene@cpht.polytechnique.fr

F. Redig  
Mathematisch Instituut  
Universiteit Leiden  
P.O. Box 9512  
2300 RA Leiden  
The Netherlands  
E-mail: redig@math.leidenuniv.nl

E. Verbitskiy  
Philips Research  
Prof. Holstlaan 4 (WO P2)  
5656 AA Eindhoven  
The Netherlands  
E-mail: evgeny.verbitskiy@philips.com